\documentstyle[12pt]{article}
\newcommand{\R}{\mbox{I$\!$R}}

\newcommand{\qed}{{\hfill {$\rlap{$\sqcap$}\sqcup$}}\\[0.2in]\hspace*{0.5in}}
\newcommand{\qedwh}{{\hfill {$\rlap{$\sqcap$}\sqcup$}}\\[0.2in]}
\newcommand{\bk}{\\[0.03in] \hspace*{0.5in} }
\newcommand{\btd}{\bigtriangledown}

\newcommand{\mfor}{\ \ \ \ {\mbox{for}} \ \ }

\begin{document}

\vspace*{0.2in}

\begin{center} {\LARGE   {\bf CONFORMAL SCALAR CURVATURE}}
\medskip \medskip \smallskip \\ {\LARGE {\bf EQUATIONS IN OPEN SPACES}} 
\medskip \medskip \medskip \\ 
{\large {Man Chun {\Large L}EUNG}}
\smallskip  \smallskip\\ 
Department of Mathematics,\\  
National University of Singapore,\\
Singapore 117543, Republic of Singapore\\
{\tt matlmc@math.nus.edu.sg}
\end{center}

\vspace{0.05in}
\begin{abstract} The article contains a brief description on the study of conformal scalar curvature equations, and discusses selected topics and questions concerning the equations in open spaces.
\end{abstract}

\vspace*{0.3in}

{\bf \Large {\bf 1. \ \ Introduction}}

\vspace*{0.2in}In nature we observe many curved objects. Apparently one of the first concepts to measure, or quantify, how curved an geometric object is, is the radius of a circle: the smaller  the radius, the more curved it is. The curvature of a curve in ${\R}^3$ is cognizant to Cauchy (cf. \cite{Willmore}).  In his generative paper ``{\it General Investigations of Curved Surfaces\,}", Gauss introduces a form of curvature, now well-known as Gaussian curvature, into the study of surfaces. The notion of curvature in $n$-dimensional spaces is propounded by Riemann in the Habilitation Lecture ``{\it On the Hypotheses Which Lie at the Foundations of Geometry\,}" and in a later work submitted to the Paris Academy in 1861. It is suggestive that a metric is not seen as inheriting from the space, but is given and can be changed without ``altering" the space. For instance, there are many metrics on $S^n$, the unit sphere in ${\R}^{n +1}$, besides the canonical metric. (For the classics of Gauss and Riemann, we refer to the book by Spivak \cite{Spivak} for informative  discussions. cf. also \cite{Chern})\\[0.05in]
\hrule

\vspace*{0.05in}
$\bullet$ KEY WORDS: conformal scalar curvature equation, positive solutions, asymptotic property.\\
$\bullet$ 2000 AMS MS CLASSIFICATIONS: Primary 35J60\,; Secondary 53C21, 58J05.\\

\newpage

\hspace*{0.5in}This critical dualism of space and metric allows for freedom in changing the metric structure, of which a simple way is to apply conformal transformations. Let $g$ be a metric on an $n$-dimensional manifold $M$, and $u$ a smooth positive function on $M$. Denote by $g_c := u^{4/(n - 2)} g$ the conformal metric for $n \ge 3$. The change in the scalar curvature, which is the second trace of the Riemannian curvature tensor, is governed by the {\it conformal scalar curvature equation} 
$$
\Delta_g u - {{n - 2}\over {4 (n - 1)}} R_g \, u + {{n - 2}\over {4 (n - 1)}} R_c \,u^{{n + 2}\over {n - 2}} = 0 \ \ \ \ {\mbox{in}} \ \ M. \leqno (1.1)
$$
Here $\Delta_g$ is the Laplacian on $(M, g)$, $R_g$ and $R_c$ are the scalar curvatures of $g$ and $g_c,$ respectively (see \cite{Besse}).\bk
The corresponding equation in dimension two takes the form
$$
\Delta_g v - K_g + K e^{2v} = 0 \ \ \ \ {\mbox{in}} \ \ M, \leqno (1.2)
$$
where $K_g$ is the Gaussian curvature of $(M, g)$, $g_c = e^{2v} g$ the conformal metric  and $K$ the Gaussian curvature of $(M, g_c)\,.$ From the works of Liouville \cite{Liouville} and Poincar\'e \cite{Poincare} to recent publications (cf. \cite{Chang-Yang.1}, \cite{Chang-Yang.2}, \cite{Chen-Ding}, \cite{Cheng-Smoller}, \cite{Kazdan-Warner}, \cite{Kiessling}, \cite{Moser}), a vast knowledge on equation (1.2) has been gained, yet a complete solution to the fundamental Nirenberg problem eludes us. The problem asks 
for sufficient and necessary conditions for a function $K$ on $S^2$ to be the Gaussian curvature function of a metric which is  conformal to the standard metric on $S^2$.\bk
The prescribed scalar curvature problem is posted and studied on compact Riemannian manifolds of higher dimension (see \cite{Aubin.2} and \cite{Kazdan}). The special case of looking for conformal metrics of {\it constant} scalar curvature is known as the Yamabe problem, which is resolved by the important efforts of Yamabe \cite{Yamabe}, Trudinger \cite{Trudinger}, Aubin \cite{Aubin.1} and Schoen \cite{Schoen.1} (see also \cite{Bahri} and \cite{Lee-Parker}). This sets a stunning example on the application of variational method and interaction between local and global geometry with analysis of partial differential equations.\bk
In a pioneer work of Loewner and Nirenberg \cite{Loewner-Nirenberg}, equation (1.1) is considered on an open domain in a compact manifold. Among other things, they show that if $\Sigma$ is a
compact $k$-dimensional submanifold of $S^n$ with $k \ge (n - 2)/2\,, \ n \ge 3\,,$ then
equation (1.1) has a smooth positive solution on $S^n \setminus \Sigma$ with $R_c$ being a negative
constant. This leads to the Singular Yamabe Problem, which asks for a complete conformal metric  of  the standard metric on $S^n \setminus \Sigma$ with constant scalar curvature, where $\Sigma$ is a compact submanifold of $S^n$. For non-positive scalar curvature, the problem is studied by 
Aviles and McOwen \cite{Aviles_McOwen}, Finn \cite{Finn}, Mazzeo \cite{Mazzeo}, and others. 
For positive scalar curvature, the problem is resolved and extended by the works of Schoen \cite{Schoen.2}, Mazzeo and Pacard (see \cite{Mazzeo-Pacard} and \cite{Pacard}; cf. also \cite{Mazzeo} and \cite{Kato}). A major tool is the subtle analysis of pseudo-differential operators on manifolds with corners developed by Melrose and others (see \cite{Melrose}; cf. also \cite{Schulze}).\bk
As we go into general open manifolds, the study becomes more difficult and diversified. Although there are other ways to see it now, the first examples to show that the Yamabe problem on complete non-compact Riemannian manifolds may not have solutions are provided by Jin \cite{Jin.1}. Despite of this, Yamabe problem may still have solutions on specific types of open manifolds (cf. \cite{Zhang}). For the case of negative scalar curvature, ample results are obtained on manifolds of negative sectional curvature in general, and hyperbolic space in particular (see, for examples,  \cite{Aviles-McOwen}, \cite{Jin.2} and \cite{Ratto-Rigoli-Veron}). In \cite{Leung.1}, \cite{Leung.2} and \cite{Leung.3}, the author studies the conformal scalar curvature equation in open Riemannian manifolds with compactifications and other general complete manifolds. \bk
Besides the complication arising from the geometry of open spaces near infinity, the 
subtlety of equation (1.1) can also be seen in one of the simplest open spaces, namely, Euclidean space ${\R}^n$. Here we confine ourselves to the case when $R_c$ is positive outside a compact subset of ${\R}^n$, $n \ge 3$ (for the case of negative conformal scalar curvature, we refer to the works of Ni and others, cf. \cite{Ni}, \cite{Li-Ni-1} and \cite{Li-Ni-2}). An useful technique called the moving plane method is introduced by Gidas, Ni and Nirenberg in \cite{Gidas-Ni-Nirenberg.1} and \cite{Gidas-Ni-Nirenberg.2}, developing an idea of Serrin on a version of the Alexandrov reflection principle. The method is found to be fruitful and deep results are obtained by Caffarelli, Gidas and Spruck \cite{Caffarelli-Gidas-Spruck}, W.-X. Chen and C.-M. Li \cite{Chen-Li}, C.-C. Chen and C.-S. Lin \cite{Chen-Lin.2} \cite{Chen-Lin.3} \cite{Chen-Lin.4}\,, and others.\bk
The apparently simple equation (1.1) calls for a variety of ideas from differential geometry,  analysis and algebra. We discuss some specific topics in the following.

\pagebreak


{\bf \Large {\bf 2. \ \ Blow-up in the Yamabe Equation}}

\vspace{0.3in}

To find a conformal metric $g_c = u^{4/(n - 2)} g$ of constant scalar curvature is equivalent to solving the {\it Yamabe equation}
$$
\Delta_gu  - {{n - 2}\over {4 (n - 1)}} R_g \, u + \lambda \,u^{{n + 2}\over {n - 2}} = 0 \ \ \ \ {\mbox{in}} \ \ M \leqno (2.1)
$$
for a constant $\lambda$, with $u$ being positive and smooth. In this section, let $(M, g)$ be a {\it compact} Riemannian $n$-manifold without boundary ($n \ge 3$)\,. We discuss a blow-up argument that is an essential part for later study. It appears in the process of the direct method in the variational problem.\bk 
For a number $s \in (0, 2n/(n - 2)]$, consider the quotient 
$$
Q_s (u) := {{ \int_M \, \left[|\btd u|^2 + {{n - 2}\over {4 (n - 1)}} R_g \,u^2 \right] \,dv_g }\over 
{ \left( \int_M |u|^s \,dv_g \right)^{2/s}   }}\,, \leqno (2.2)
$$
where $u \in L^2_1 (M) \setminus \{ 0 \}\,.$  Let 
$$
\lambda_s = \inf \, \{ Q_s (u)  \ | \ u \in L^2_1 (M) \setminus \{ 0 \} \} \ \ \ \ {\mbox{and}} \ \ \lambda (M) = \lambda_{{2n}\over{n - 2}}\,.  
$$
By using a direct minimizing procedure, it can be shown that for $2 < s < 2n/(n - 2)$, there exists a smooth positive function $u_s$ such that its $L^s$-norm is equal to one, $Q_s (u_s) = \lambda_s$ and $u_s$ satisfies the equation
$$
\Delta_gu_s  - {{n - 2}\over {4 (n - 1)}} R_g \, u_s + \lambda_s \,u_s^{s - 1} = 0 \ \ \ \ {\mbox{in}} \ \ M\,. \leqno (2.3)
$$
The direct method does not work when $s = 2n/(n - 2)$ because the Sobolev embedding $L^2_1 (M) \hookrightarrow L^{2n/(n - 2)} (M)$ is continuous but not compact. However, if one can show that there exists a positive constant $c$ such that $u_s \le c$ in $M$ for $2 < s < 2n/(n - 2)$, then there exists a sequence $\{ s_i \} \in (2, 2n/(n - 2))$ such that $s_i \to 2n/(n - 2)$ and 
$u_{s_i}$ converges to a smooth positive function $u$ which satisfies the Yamabe equation  (2.1).\bk
Suppose that no such upper bound $c$ exists. It follows that there exist sequences $\{ s_k \} \subset (2, \,2n/(n - 2))$ and $\{ y_k \} \subset M$ such that 
$$
s_k \to {{2n}\over {n - 2}} \ \ \ \ {\mbox{and}} \ \ \ \ m_k := u_{s_k} (y_k) = \max u_{s_k} \to \infty
 \ \ \ \ {\mbox{as}} \ \ \ \ k \to \infty.
$$ 
As $M$ is compact, we may assume that $y_k \to y_o$ as $k \to \infty$. For a normal coordinate system centered at $y_o$ and with radius $\rho$, let the coordinates of $y_k$ be $x_k$, $k = 1, 2,....$ In the local coordinates, $u_{s_k}$ satisfies the equation 
$$
{1\over { \sqrt{ {\mbox{det}} \, g } }}\, \partial_j \left(\sqrt{ {\mbox{det}} \,g } \, g^{ij} \,\partial_i u_{s_k} \right)  - {{n - 2}\over {4 (n - 1)}} R_g \, u_{s_k} + \lambda_{s_k} \,u_{s_s}^{s_k - 1} = 0 \ \ \ \ {\mbox{in}} \ \ B_o (\rho)\,. \leqno (2.4)
$$
\hspace*{0.5in}The idea here is to consider the normalized function
$$
v_k (x) := {{u_k (\delta_k x + x_k)}\over {m_k}}\,, \ \ \ \ {\mbox{where}} \ \ \delta_k := {1\over {m_k^{(s_k - 2)/2}}}\,. 
$$
We have $x_k \to 0$ and $\delta_k \to 0$ as $k \to \infty$. Here $v_k$ is defined on a ball in ${\R}^n$ of radius $\rho_k = (\rho - |x_k|)/\delta_k$ and $\rho_k \to \infty$ as $k \to \infty$.  
By the argument of diagonal subsequence and the property of normal coordinates (see \cite{Schoen-Yau-Book}), one observes that a subsequence of $\{ v_k \}$ converges to a smooth positive function $v$ which is a solution of the equation 
$$
\Delta_o\, v + \lambda \,v^{{n + 2}\over {n - 2}} = 0 \ \ \ \ {\mbox{in}} \ \ {\R}^n, \leqno (2.5)
$$
where $\lambda = {\mathop{\lim}\limits_{k \to \infty}} \lambda_{s_k}$ and  $\Delta_o$ is the standard Laplacian on ${\R}^n$. It is known that $\lambda \le \lambda (M)$ if $\lambda (M) < 0$; and $\lambda = \lambda (M)$ if $\lambda (M) \ge 0$ \cite{Schoen-Yau-Book}. Let $d$ be the diameter of $(M, g)$. By a change of variables we have
$$
\int_{|x| \le {d\over 2} \delta_k^{-1} } v_k^{s_k} \sqrt{ {\mbox{det}} \, g} \, dx =  \delta_k^{ {{ 2s_k}\over {s_k - 2}} - n} \int_{B_{x_k} ({d\over 2})} u^{s_k}_k \, dv_g  \le \delta_k^{ {{ 2s_k}\over {s_k - 2}} - n}\,, \leqno (2.6)
$$
where we note that 
$$
{{ 2s_k}\over {s_k - 2}} - n > 0 \ \ \ \ {\mbox{and}} \ \ \to 0 \ \ \ \ {\mbox{as}} \ \ k \to \infty\,.
$$
Here $B_{x_k} ({d/2})$ denotes the open ball in $(M, g)$ with center at $x_k$ and radius equal to $d/2$. 
From (2.6), the Fatou lemma and $\displaystyle{\lim_{k \to \infty} \delta_k = 0}\,,$  
we obtain 
$$
\int_{\R^n} v^{{2n}\over {n - 2}} \, dx \le 1\,.  \leqno (2.7)
$$
We see that if $\lambda (M) \le 0$, then $v$ is a positive subharmonic function on ${\R}^n$ and $v \in L^{2n/(n - 2)} ({\R}^n).$ Hence $v\equiv 0$ (a Liouville-type theorem, see \cite{Schoen-Yau-Book}), which is a contradiction. The conclusion can be drawn that the upper bound $c$ exists if $\lambda (M) \le 0$. We consider the case $\lambda (M) > 0$ in the next section.

\pagebreak


{\bf \Large {\bf 3. \ \ Global Solutions}}

\vspace{0.3in}

In case $\lambda (M) > 0$, by the scaling 
$$
u = \left[ {{\lambda (M)}\over {n (n - 2)}} \right]^{{n - 2}\over 4} v\,, \leqno (3.1)
$$
equation (2.5) can be rewritten as 
$$
\Delta_o\, u + n (n - 2) \,u^{{n + 2}\over {n - 2}} = 0 \ \ \ \ {\mbox{in}} \ \ {\R}^n. \leqno (3.2)
$$
An amazing result of  Gidas, Ni and Nirenberg \cite{Gidas-Ni-Nirenberg.1} \cite{Gidas-Ni-Nirenberg.2} exhibits all smooth positive global solutions of equation (3.2).\\[0.2in]
{\bf Theorem 3.3.} \ \ {\it Let $u$ be a smooth positive solution of equation (3.1) in ${\R}^n, \ n \ge 3$. 
Then $u$ is rotationally symmetric about a point $\xi \in {\R}^n\,.$ Furthermore,  
$$
u (x) = \left( {{\lambda}\over {\lambda^2 + |x - \xi|^2}} \right)^{{n  - 2}\over 2}\ \ \ \ for \ \ x \in {\R}^n \leqno (3.4)
$$
and for a positive constant $\lambda$.}

\vspace*{0.2in}

\hspace*{0.5in}We observe that if $u$ is given by (3.4), then the conformal transformation $u^{4/(n - 2)} g_o$ changes the Euclidean metric $g_o$ into a spherical metric on $S^n$. Hence solutions in (3.4) are known as {\it spherical solutions}, and are denoted by $u_{\lambda, \,\xi}$ for $\lambda > 0$ and $\xi \in {\R}^n$.\bk
If we couple theorem 3.3 with the blow-up process we may draw some more conclusions. 
A calculation shows that 
$$
\int_{\R^n} u_{\lambda, \,\xi}^{{2n}\over {n - 2}} \, dx = {{\omega_n}\over {2^n }}  \leqno (3.5)
$$
for $\lambda > 0$ and $\xi \in {\R}^n$. Here $\omega_n$ is the volume of $S^n$ in ${\R}^{n + 1}$ with the standard metric. (The notation may be different from some other works.) It follows from (2.7), (3.1) and (3.5) that 
$$
1 \ge \int_{\R^n} v^{{2n}\over {n - 2}} \, dx = \left[ {{n (n - 2)}\over {4 \lambda (M)}} \right]^{n\over 2} \omega_n\,.
$$
Hence we obtain
$$
\lambda (M) \ge {{n (n - 2)}\over 4}\, \omega_n^{2\over n}.
$$
The estimate is not sharp though. Indeed, by using the Sobolev inequality and equations (2.5) and (3.5), we obtain
$$
\Lambda \left( \int_{{\R}^n} v^{{2n}\over {n - 2}} \, dx \right)^{{n - 2}\over n} \le 
\int_{{\R}^n} |\btd v|^2 \, dx = \lambda (M) \int_{{\R}^n} v^{{2n}\over {n - 2}} \, dx\,. \leqno (3.6)
$$
It follows from (2.7) that $\lambda (M) \ge \Lambda$, and the best constant $\Lambda$ in the Sobolev inequality is found to be equal to $\lambda (S^n) = n (n - 1)\, \omega_n^{2\over n}$. Therefore, under 
the assumption that the uniform upper bound $c$ does not exist, we have 
$$
\lambda (M) \ge n (n - 1)\, \omega_n^{2\over n}\,.
$$
We are led to the conclusion that the Yamabe problem has a solution if $\lambda (M) < \lambda (S^n)$. 
The major remaining part is to show that $\lambda (M) < \lambda (S^n)$ {\it unless} $(M, g)$ is conformally equivalent to $S^n$ (see \cite{Schoen-Yau-Book} and \cite{Lee-Parker}). Notice that the exact form of $u$ is not utilized in (3.6)


\vspace*{0.4in}

{\bf \Large {\bf 4. \ \ Non-constant Positive Scalar Curvature}}

\vspace{0.3in}

Consider the conformal scalar curvature equation
$$
\Delta  u + n (n - 2) K u^{{n + 2}\over {n - 2}} = 0 \ \ \ \ {\mbox{in}} \ \ {\R}^n, \leqno (4.1)
$$
where $K$ is a smooth function on ${\R}^n$ with  
$$
a^2 \le K (x) \le b^2 \mfor |x| \gg 1\,. \leqno (4.2)
$$
Here $a$ and $b$ are positive constants. Of fundamental concerns are the asymptotic behaviors of $u (x)$ when $|x| \to \infty$. The basic asymptotic properties we seek to consider are the following.\\[0.1in]
{\bf (i) \ \ Fast decay:} \ \ $u (x) \le C |x|^{- (n - 2)}$ for $|x| \gg 1\,.$ \\[0.05in] 
{\bf (ii) \ \ Slow decay:} \ \ $u (x) \le C |x|^{- (n - 2)/2}$ for $|x| \gg 1\,.$ \\[0.05in]
{\bf (iii) \ \ Lower bound/``Completeness":} \ \ $u (x) \ge C |x|^{- (n - 2)/2}$ for $|x| \gg 1\,.$\\[0.1in]
\hspace*{0.5in}We first observe that the bound $u (x) \ge C |x|^{- (n - 2)}$ for $|x| \gg 1$ is guaranteed by the maximum principle. Denote by $B_o (r)$ the open ball in ${\R}^n$ with radius $r > 0$ and center at the origin.\\[0.2in]
{\bf Lemma 4.3.} \ \ {\it Let $u$ be a smooth positive superharmonic function (i.e., $\Delta u \le 0$) on ${\R}^n \setminus \overline{B_o (r_o)}\,.$ There exist positive numbers $c$ and $r_1$ such that} 
$$
u (x) \ge c |x|^{2 - n} \ \ \ \ for \ \ |x| \ge r_1\,.
$$

\vspace*{0.1in}

{\bf Proof.} \ \ Fix a number $\rho > r_o$. We can find a positive number $c$  such that 
$$
u (x) \,\ge \,c \,|x|^{2 - n}  \ \ \ \ \mfor \  \ |x| = \rho\,.
$$
Let 
$$
\phi (x) =  c |x|^{2 - n} \ \ \ \ {\mbox{and}} \ \ w (x) = u (x) - \phi (x) \mfor x \in {\R}^n \setminus \overline{B_o (r)}\,.
$$
Suppose that there is a point $x_o$ with $|x_o| > \rho$ such that $w (x_o) < 0\,,$ i.e., $\phi (x_o) > u (x_o)\,.$ Let $\delta = -w (x_o)\,.$ There exists a number $R > |x_o|$ such that 
$$
\phi (x) \le {\delta \over 2} \mfor |x| = R\,.
$$
Hence 
$$
w (x) \ge 0 \mfor |x| = \rho \ \ \ \ {\mbox{and}} \ \ w (x) > - {\delta\over 2} \mfor |x| = R\,.
$$
Since $w (x_o) = - \delta\,$,\,   
inf $\{ w \ | \ B_o (R) \setminus \overline{B_o (\rho)}\,\}$ is achieved at a point inside $B_o (R) \setminus \overline{B_o (\rho)}\,.$ As
$$
\Delta w = \Delta u \le 0 \ \ \ \ {\mbox{in}} \ \ 
{\R}^n \setminus \overline{B_o (r)}\,,
$$
we are led to a contradiction via the maximum principle. Hence $w (x) \ge 0$ for all $|x| \ge \rho\,.$\qed
The lemma is sharp as it can be seen that the spherical solutions decay in the order $O (|x|^{2-n})$ for $|x| \gg 1$. The analytic consequences of the slow decay can be seen through the following lemma, which is observed in [18] (a version of the lemma is proved in [8]). 
\\[0.2in]
{\bf Lemma 4.4.} \ \ {\it Assume
that
$K$ satisfies (4.2). 
 Let $u$ be a smooth positive solution of (4.1)
in ${\R}^n$. Then the following statements are equivalent.\\[0.05in]
{\bf (A)} \ \ $u$ has slow decay (ii).\\[0.05in]
{\bf (B)} \ \ $u$ satisfies a spherical Harnack inequality, i.e., there exists
a positive constant $C_h$  such that the inequality 
$$
\sup_{S_r} u \, \le  \,C_h\, \inf_{S_r} u   $$
holds for all $r > 0\,.$ Here $S_r$ is the sphere in ${\R}^n$ with center at the origin and radius equal to $r$. \\[0.05in]
{\bf (C)} \ \ There exists a positive constant $C_g$  such that $|x| \cdot
|\btd u (x)| \le C_g\, u (x)$ for $|x| > 0\,.$}\\[0.2in]
\hspace*{0.5in}By using the {\it Kelvin transformation}, we can bring, as it were, the ``infinity" to the origin. Let
$$
y = {{x}\over {|x|^2}} \ \ \ \ {\mbox{and}} \ \ \ \   
w (y) := |y|^{2 - n} \,u \left( {y\over {|y|^2}} \right) 
\mfor x, \, y \in {\R}^n \setminus \{ 0 \}\,. \leqno (4.5)
$$
Then $w$ satisfies the equation 
$$
\Delta w (y) + n (n - 2) \,K_o (y)\, w^{{n + 2}\over {n - 2}} (y) = 0 \mfor y \not= 0\,, \leqno (4.6)
$$
where $K_o (y) = K (y/|y|^2)$ for $y \not= 0\,.$ In this way, the questions on the asymptotic properties of $u$ may be paraphrased as questions on the singularity/regularity of $w$ at the origin. In particular,\\[0.1in]
{\bf (a)} \ \ $w$ has a removable singularity at $0$ 
\ \ $\Longrightarrow$ \ \  $u$ has {\bf fast decay (i)}.\\[0.05in] 
{\bf (b)} \ \ $w (y) \le C |y|^{- (n - 2)/2}\,$ for $|y| \sim 0$, $|y| \not= 0$ \ \ $\Longleftrightarrow$ \ \ 
$u$ has {\bf slow decay (ii)}.\\[0.05in] 
{\bf (c)} \ \ $w (y) \ge C |y|^{- (n - 2)/2}\,$ for $|y| \sim 0$, $|y| \not= 0$ \ \ $\Longleftrightarrow$ \ \ $u$ has {\bf lower bound (iii)}.

\vspace*{0.4in}

{\bf \Large {\bf 5. \ \ Slow Decay and Blow-up}}

\vspace{0.3in}

We examine the slow decay and blow-up of smooth positive solutions $w$ of the equation
$$
\Delta w + n (n - 2) K_o \,w^{{n + 2}\over {n - 2}} = 0 \ \ \ \ {\mbox{in}} \ \ B_o (1) \setminus \{ 0 \}\,. \leqno (5.1)
$$
Given a small positive number $\varepsilon > 0\,,$ define
$$
d_\varepsilon (x) = \max \, \left\{ 0, \ \min\, \{ |x| - \varepsilon\,, \ 5/8 - |x| \} \right\} 
\mfor x \in B_o (5/8)\,.
$$
Set 
$$
W_\varepsilon (x) =  [d_\varepsilon (x)]^{{n - 2}\over 2}  \,w (x) \mfor x \in B_o (5/8)\,, \leqno (5.2)
$$
and 
$$
M_\varepsilon = \sup_{B_o (5/8)} W_\varepsilon\,.
$$
It follows that if there exists a positive constant $C$ such that for all $x \in B_o (5/8)$, 
$$
M_\varepsilon (x) \le C \ \ \ \ {\mbox{for \ all}} \ \ \varepsilon > 0 \ \ {\mbox{small}}\,, \leqno (5.3)
$$
($C$ is independent on $x$ and $\varepsilon$), then $w$ has slow decay in the form of (b). Assume that no such constant $C$ exists, i.e., there is a decreasing sequence of small positive numbers $\{ \varepsilon_i \}$ such that $M_{\varepsilon_i} \to \infty$ as $i \to \infty\,.$  Let $x_i$ be a maximum point of $W_{\varepsilon_i}\,,$ $i = 1, \,2,....$  Suppose that $|x_i| > c$ for a positive number $c$ and for $i = 1, \,2,....$ We have  
$$
M_{\varepsilon_i} \le (5/8)^{{n - 2}\over 2} \sup_{B_o (5/8) \setminus B_o (c)} \,w \ \ \ \ {\mbox{for}} \ \ i = 1, \,2,...,  \leqno (5.4)
$$
where the right hand side is a bounded number. But this is a contradiction. Therefore, by choosing a subsequence if necessary, we may assume that ${\mathop{\lim}\limits_{i \to \infty}} x_i = 0\,.$ Similar to the blow-up procedure in section 2 (see \cite{K-M-P-S}), let 
$$
\mu_i = {1\over {[w (x_i)]^{2/(n - 2)} }}
$$
and 
$$
v_i \,(x) = \mu_i^{{n - 2}\over 2} w (x_i + \mu_i \,x)\,, \ \ \ \ i = 1, \,2,.... \leqno (5.5)
$$
We assert that, under the conditions
$$
\lim_{y \to 0} K_o (y) = 1 \ \ \ \ {\mbox{and}} \ \ \ \ |y| \cdot |\btd_y K_o (y)| \le C_o \mfor y \in B_o (5/8)\,, \leqno (5.6)
$$
there is a convergent subsequence of $\{ v_i \}$. Here $C_o$ is a constant. To proceed, set 
$$
r_i = d_{\varepsilon_i} (x_i)/2\,, \ i = 1, \,2,....
$$ 
For $y \in B_{x_i} (r_i)\,,$ we have 
$$
W_{\varepsilon_i} (y) \le W_{\varepsilon_i} (x_i) \ \ \ \ {\mbox{and}} \ \ d_{\varepsilon_i} (y)
 \ge {{ d_{\varepsilon_i} (x_i)}\over 2}\,, \ \ \ \ i = 1, \,2,....
$$
It follows that $w (y) \le 2^{(n - 2)/2} w (x_i)$ for $y \in B_{x_i} (r_i)$. That is, 
$$
w (x_i + \mu_i \,x) \le 2^{(n - 2)/2} w (x_i) \mfor x \in B_o (r_i/\mu_i)\,. \leqno (5.7)
$$ 
Let 
$$
R_i = {{r_i}\over {\mu_i}} = {1\over 2} M_i^{2/(n - 2)}.
$$
From (5.7), we have
$$
v_i (x) \le 2^{(n - 2)/2} \mfor x \in B_o (R_i) \ \ {\mbox{and}} \ \ i = 1, \,2,....
$$
Furthermore, $v_i$ satisfies the equation
$$
\Delta_x v_i (x) + n (n - 2) K_o\, (x_i + \mu_i \,x)\, v_i^{{n + 2}\over {n - 2}} (x) = 0 \ \ \ \ {\mbox{in}} \ \ B_o (R_i)\,, \ i = 1, \,2,....
$$
We have
$$
|\btd_x K_o\, (x_i + \mu_i \,x)| = \mu_i |\btd_y K_o\, (y)|\vert_{y = x_i + \mu_i \,x}  \mfor x \in B_o (R_i)\,.
$$
Observe that 
\begin{eqnarray*}
& \, & \mu_i |\btd_y K_o \,(y)|\vert_{y = x_i + \mu_i \,x}\\
 & = & {{|\btd_y K_o\, (y)|\vert_{y = x_i + \mu_i \,x} }\over {[u (x_i)]^{{n - 2}\over 2} }} \le d_{\varepsilon_i} (x_i) |\btd_y 
K_o\, (y)|\vert_{y = x_i + \mu_i \,x}\\
& \le & 2\, d_{\varepsilon_i} (y) |\btd_y K_o\, (y)|\vert_{y = x_i + \mu_i \,x} \le 2 |y| \cdot  |\btd_y K_o\, (y)|\vert_{y = x_i + \mu_i \,x}
\end{eqnarray*}
for $y \in B_{x_i} (r_i)$ and large $i$.  
As $i \to \infty\,,$ $x_i \to 0\,,$ $\mu_i \to 0\,,$ $\mu_i R_i = r_i \to 0$ and $R_i \to \infty$, thus if we assume (5.6), 
then standard elliptic theory implies that there is a subsequence of $\{ v_i \}$ which converges in $C^2$ norm on compact subsets to a $C^2$ function $u$ which satisfies the equation
$$
\Delta u + n (n - 2)\,u^{{n + 2}\over {n - 2}} = 0 \ \ \ \ {\mbox{in}} \ \ {\R}^n. \leqno (5.8) 
$$
If it {\it were} the case that the above equation had no positive solutions, then $w$  would have had slow decay in the form of (b). But, as we see, equation (5.8) have a family of positive solutions given by (3.4). 

\pagebreak


{\bf \Large {\bf 6. \ \ Reflections}}

\vspace{0.3in}

We discuss the reflection argument leading to slow decay when $K_o \equiv 1$ in a neighborhood of the origin. Suppose  that  the constant $C$ in (5.3) does not exist. In the blow-up analysis of equation (5.1),  given any positive number $\delta$, we can find a $\varepsilon > 0$ and a maximal point $x_1$ of $W_\varepsilon$ such that 
$$
\| v_\mu - u_{1, o} \|_{C^2 (B_o (R))} < \delta  \leqno (6.1)
$$
(see \cite{K-M-P-S}). Here $u_{1, o}$ is the spherical solution given by (3.4) with $\lambda = 1$ and $\xi$ being the origin, and  
$$
v_\mu (x) := \mu^{{n - 2}\over 2} w (x_1 + \lambda x) \ \ \ \ {\mbox{with}} \ \ \mu = [w (x_1)]^{2 \over {2 -n}}. \leqno (6.2)
$$
The reflection process appears to be geometrically natural in the {\it extended cylindrical coordinates}. Let 
$$
t = - \ln |x| \ \ \ \  {\mbox{and}} \ \ \ \theta = x/|x| \mfor x \not=0\,. \leqno (6.3)
$$ 
Define
$$
v (t, \theta) = |x|^{{n - 2}\over 2} v_\mu (x) \leqno (6.4)
$$
with $x$ represented by $(e^{-t}, \ \theta)$ in the polar coordinates. For not too negative $t$,  $v$ satisfies the equation 
$$
{{\partial^2 v}\over {\partial s^2}} + \Delta_\theta\, v - \left( {{n -
2}\over 2} \right)^2 v + n (n - 2) v^{{n + 2}\over {n - 2}} = 0\,, \leqno (6.5)
$$ 
where $\Delta_\theta$ is the standard Laplacian on $S^{n - 1}$. 
Define the differential operator $L$ by  
$$
L = {{\partial^2 }\over {\partial s^2}} + \Delta_\theta\,  - \left( {{n -
2}\over 2} \right)^2\,. \leqno (6.6)
$$
The following lemma can be derived from the mean value inequality and the Hopf lemma (cf. \cite{Gidas-Ni-Nirenberg.1} and \cite{K-M-P-S}).\\[0.2in]
{\bf Lemma 6.7.} \ \ {\it Let $u$ and $u^*$ be $C^2$ functions which satisfy the equations
$$
L u + n (n - 2) K u^{{n + 2}\over {n - 2}} = 0 \ \ \ \ {\mbox{and}} \ \ \ \ 
L u^* + n (n - 2) K^* {u^*}^{{n + 2}\over {n - 2}} = 0\,,
$$
respectively, in an open connected set $\Omega \subset R \times S^{n - 1}.$ Here $K$ and $K^*$ are smooth functions on $\Omega$ with $K^* \le K$ in $\Omega$. Assume that $u^* \le u$ in $\Omega$. Then we have the following conclusions.}

\begin{enumerate}

\item {\it Either $u^* < u$ or $u^* \equiv u$ in $\Omega$.}
\item {\it Suppose that there is a ball $B$ in $\Omega$ and a point $p \in \partial B \cap \partial \Omega$ with $u (p) = u^* (p)$. Assume also that $u$ and $u^*$ are continuous in $\Omega \cup \{ p \}\,.$ Then if $u \not\equiv u^*$ in $B$, we have} 
$$
{{\partial (u - u^*)}\over {\partial {\bf n} }}  < 0\,,
$$
{\it where ${\bf n}$ is the unit outward normal of $B$ at $p$.}

\end{enumerate}
\vspace*{0.02in}
\hspace*{0.5in}From (6.1) we see that $v$ is close in $C^2$ norm on $[- \ln R\,, \ \infty) \times S^{n - 1}$ to the function $v_o (t, \theta) := (2\cosh t)^{(2- n)/2}.$ One finds that 
\begin{eqnarray*}
{{\partial v_o (t, \theta)}\over {\partial t}}  & = & \left( {{n - 2}\over 2} \right) {{e^{-t} - e^t}\over { (2\cosh t)^{n\over 2} }} \\
{\mbox{and}} \ \ \ \ 
{{\partial^2 v_o (t, \theta)}\over {\partial t^2}}  & = &  ( n - 2 )\, 
{{n \sinh^2 t - 2 \cosh^2 t}\over { (2\cosh t)^{{n\over 2} +1} }}\,.  
\end{eqnarray*}
It follows that there exists positive numbers $\varepsilon$ and $\delta$ such that for any $\theta \in S^{n - 1}$ we have
$$
{{\partial v_o (\varepsilon, \theta)}\over {\partial t}}  < 0\,, \ \ \ \ {{\partial v_o (-\varepsilon, \theta)}\over {\partial t}}  > 0 \ \ \ \ {\mbox{and}} \ \ \ \ {{\partial^2 v_o (t, \theta)}\over {\partial t^2}} < \delta \mfor t \in [-\varepsilon, \varepsilon]\,.
$$
As $v$ can be assumed to be arbitrarily close in $C^2$ norm on $[- \ln R\,, \ \infty) \times S^{n - 1}$ to $v_o$, it follows that, given $\theta \in S^{n - 1}$, there is $t_o \in (-1, 1)$ such that 
$$
{{\partial v(t_o, \theta)}\over {\partial t}} = 0 \ \ \ {\mbox{and}} \ \ {{\partial^2 v (t_o, \theta)}\over {\partial t^2}} < 0\,. \leqno (6.8)
$$
For a number $t_1$, we define the reflection across the sphere $\{ t_1 \} \times S^{n - 1}$ by 
$$
v^* (t, \,\theta) = v (2 t_1 - t, \,\theta)\,. \leqno (6.9)
$$
Based on the graph of the function $v_o$, for $t_1$ very large, we have 
\begin{eqnarray*}
\inf\, \{ v ( - \ln \left( {1\over {16 \lambda}} \right)\,, \theta)  \ | 
\ \theta \in S^{n - 1} \} & > & 
\sup \,\{ v (2t_1 + \ln \left( {1\over {16 \lambda}} \right)\,, \, \theta)  \ | \ \theta \in S^{n - 1} \}\\
& = & \sup \,\{ v^* ( - \ln \left( {1\over {16 \lambda}} \right)\,, \,\theta)  \ | \ \theta \in S^{n - 1} \}\,.
\end{eqnarray*}
An reflection argument originally used by Alexandrov (cf. Ch. 7 in \cite{Hopf}) provides useful relations between the ``righthand side" and the ``lefthand side", especially at a value $t_o \in (-1, 1)$ which satisfies (6.8) (cf. lemma 6.7). It is asserted in \cite{K-M-P-S} that 
$$
\inf \,\{ v ( - \ln \left( {1\over {16 \lambda}} \right)\,, \theta)  \ | \ \theta \in S^{n - 1} \} < 
\sup \,\{ v (t \,, \theta)  \ | \ \ t > \ln \left( {1\over {16 \lambda}} \right) - 2\,,\ \theta \in S^{n - 1} \}\,.
$$
It can be seen from the definitions of $\lambda$ and $v$ that 
$$
\inf \{ v ( - \ln \left( {1\over {16 \lambda}} \right)\,, \,\theta)  \ | \ \theta \in S^{n - 1} \}= 16^{{2 - n}\over 2} \inf_{\partial B_{x_1} (1/16)} w\,.
$$
As $v$ is close to $v_o$ for $t > 0$, and $w$ is a positive superharmonic function in $B_o (1)\setminus \{ 0\}$ with a singularity at $0$, we obtain 
$$
\inf_{\partial B_o (3/4)} w \le \inf_{\partial B_{x_1} (1/16)} w \le C_o v_o (\ln \left( {1\over {16 \lambda}} \right) - 2) \le C \lambda^{{n - 2}\over 2} = C w^{-1} (x_1)\,,
$$
where $C_o$ and $C$ are constants. That is,
$$
w (x_1) \le C \,[\inf_{\partial B_o (3/4)} w]^{-1}. \leqno (6.10)
$$
But this contradicts that $u (x_1)$ can be chosen to be arbitrarily large. 
It follows that when $K_o\equiv 1$, smooth positive solutions to equation (4.1) have slow decay (b).

\vspace*{0.4in}

{\bf \Large {\bf 7. \ \ Spherical Solutions}}

\vspace{0.3in}

Spherical solutions appear to be essential in the study of equation (4.1). We consider some of their basic properties. Given spherical solutions 
$$
u_1 (x) = \left( {\lambda_1\over { \lambda_1^2 + |x - \xi_1|^2}} \right)^{{n - 2}\over 2}
\ \ \ \ {\mbox{and}} \ \ \ \ 
u_2 (x) = \left( {\lambda_2\over { \lambda_2^2 + |x - \xi_2|^2}} \right)^{{n - 2}\over 2} \leqno (7.1)  
$$
for $x \in {\R}^n$, let 
$$
v_1 (s, \theta) = |x|^{{n - 2}\over 2} u_1 (x) \ \ \ \ {\mbox{and}}  \ \ \ \ v_2 (s, \theta) = 
|x|^{{n - 2}\over 2} u_2 (x)\,, \leqno (7.2) 
$$
where $e^s = |x|$ and $\theta = x/|x|$ for $x \in {\R}^n \setminus \{ 0 \}\,.$ We call $v_1$ and $v_2$ the {\it cylindrical transforms} of $u_1$ and $u_2$, respectively. Note that there is a sign change comparing to (6.3).\\[0.2in]
{\bf Proposition 7.3.} \ \ {\it In (7.1), if  $\xi_1/\lambda_1 = \xi_2/\lambda_2$, then}
$$
v_1 (s, \theta) = v_2 (s + \bar s, \theta) \ \ \ \ for \ \ s \in \R \ \ and \ \ \theta \in S^{n - 1}, \ \ where \ \ \bar s = \ln\, (\lambda_2/\lambda_1)\,.
$$

\vspace*{0.1in}

{\bf Proof.} \ \  We have

\pagebreak

\begin{eqnarray*}
& \ & v_1 (s, \theta)\\ 
&=& e^{{{n - 2}\over 2} s}  \left( {\lambda_1\over { \lambda_1^2 + |x - \xi_1|^2}} \right)^{{n - 2}\over 2}\\
 & = & \left( {{\lambda_1 e^s}\over {\lambda_1^2 + |x|^2 - 2 x \cdot \xi_1 + |\xi_1|^2}} \right)^{{n - 2}\over 2}
= \left( {{\lambda_1 e^s}\over {\lambda_1^2 + e^{2s} - 2 \, e^s (\theta \cdot \xi_1) + |\xi_1|^2}} \right)^{{n - 2}\over 2}\\
& = & \left( {{ {{e^s}\over {\lambda_1}} }\over {1 + {{e^{2s}}\over {\lambda_1^2}} 
- 2 {{e^s}\over \lambda_1}\, \left( \theta \cdot {\xi_1\over \lambda_1} \right) + {{|\xi_1|^2}\over {\lambda_1^2}} }} \right)^{{n - 2}\over 2}
=\left( {{ {{e^s}\over {\lambda_1}} }\over {1 + {{e^{2s}}\over {\lambda_1^2}} 
- 2 {e^s\over \lambda_1}\, \left( \theta \cdot {\xi_2\over \lambda_2} \right) + {{|\xi_2|^2}\over {\lambda_2^2}} }} \right)^{{n - 2}\over 2}\\
& = & \left( {{ \lambda_2\, \left( {{\lambda_2 \,e^s}\over {\lambda_1}} \right)}\over {\lambda_2^2 + {{\lambda_2^2\,e^{2s}}\over {\lambda_1^2}} 
- 2 {\lambda_2 \,e^s\over \lambda_1}\, \left( \theta \cdot \xi_2 \right) + |\xi_2|^2 }} \right)^{{n - 2}\over 2}\\
& = &\left( {{ \lambda_2 \, e^{s + \bar s}}\over {\lambda_2^2 +  e^{2(s + \bar s)} 
- 2 \,e^{s + \bar s}\, \left( \theta \cdot \xi_2 \right) + |\xi_2|^2 }} \right)^{{n - 2}\over 2}\,.
\end{eqnarray*}
On the other hand, 
$$
v_2 (s\,, \theta) = \left( {{\lambda_2 e^s}\over {\lambda_2^2 + e^{2s} - 2 \, e^s (\theta \cdot \xi_2) + |\xi_2|^2}} \right)^{{n - 2}\over 2}\,.
$$
Hence we obtain the result. (Here $\cdot$ denotes the dot product for vectors in ${\R}^n$.) \qed
One may think of spherical solutions as represented by points in the upper half space
$$
{\R}^{n + 1}_+ = \{ (\xi_1,..., \xi_n, \lambda) \ | \ \xi = (\xi_1,..., \xi_n) \in {\R}^n, \ \lambda \in  {\R}^+ \}\,.
$$
Introducing the relation $\xi/\lambda = \xi'/\lambda'$, we obtain the upper hemisphere $S_+^{n + 1}\,.$ The north pole corresponds to the class with $\xi = 0$. \bk
For a function $u$ defined on ${\R}^n$, the {\it Kelvin transform about the sphere of radius $a > 0$} is given by 
$$
{\tilde u} (x) = |x|^{2 - n} \, a^{n - 2} \, u \left( {{a^2 x}\over {|x|^2}} \right) \mfor x \in {\R}^n \setminus \{ 0 \}\,. \leqno (7.4)
$$
In case $u$ is a smooth positive function on ${\R}^n$ which satisfies equation (4.1), a direct calculation shows that $\tilde u$ satisfies the equation
$$
\Delta \tilde u (x) + n (n - 2) K (a^2x/|x|^2)\, {\tilde u}^{{n + 2}\over {n - 2}} (x) = 0 \mfor x \in {\R}^n \setminus \{ 0 \}\,. \leqno (7.5)
$$
If $K \equiv 1$ in ${\R}^n$ and $u$ is a spherical solution, then, from (7.5) and the observation that $\tilde u$ has a removable singularity at $0$,  we see via theorem 3.3 that $\tilde u$ is also a spherical solution.\\[0.2in]
{\bf Proposition 7.6.} \ \ {\it Let $u$ be a spherical solution given by}
$$
u (x) = \left( {\lambda\over { \lambda^2 + |x - \xi|^2}} \right)^{{n - 2}\over 2} \ \ \ \ for \ \ x \in {\R}^n,
$$
{\it and $\tilde u$ the Kelvin transform of $u$ about the sphere of radius $a > 0$. We have}
$$
\tilde u = \left( {\tilde \lambda\over { {\tilde  \lambda}^2 + |x - \tilde  \xi|^2}} \right)^{{n - 2}\over 2} \ \ \ \ for \ \ x \in {\R}^n, 
$$
{\it where}
$$
\tilde \lambda = {{a^2 \lambda}\over {\lambda^2 + |\xi|^2}} \ \ \ \ and \ \ \ \ \tilde \xi = 
{{a^2 \,\xi}\over {\lambda^2 + |\xi|^2}}\,. \leqno (7.7)
$$
{\it In particular, we have $\tilde \xi /\tilde \lambda = \xi/\lambda\,.$}\\[0.1in]
{\bf Proof.} \ \ From (7.4) we obtain
\begin{eqnarray*}
\tilde u (x) & = & |x|^{2 - n} \, a^{n - 2} \, 
\left( {\lambda\over { \lambda^2 + \bigg\vert {{a^2 x}\over{ |x|^2}}- \xi \bigg\vert^2 }} \right)^{{n - 2}\over 2}\\
& = & \left( {{ a^2 \lambda}\over { |x|^2 \lambda^2 +  a^4 - 2 a^2 \xi \cdot x +  |\xi|^2 |x|^2}} \right)^{{n - 2}\over 2}\\
& = & \left[ {{ a^2 \lambda/(\lambda^2 + |\xi|^2)}\over { {{a^4}\over {\lambda^2 + |\xi|^2}} + |x|^2  - {{2 a^2 \xi \cdot x}\over {\lambda^2 + |\xi|^2}} }} \right]^{{n - 2}\over 2}\\
& = & \left[ {{ a^2 \lambda/(\lambda^2 + |\xi|^2)}\over { {{a^4 \lambda^2}\over {(\lambda^2 + |\xi|^2)^2}} + \left( x - {{a^2 \,\xi}\over {\lambda^2 + |\xi|^2 }}\right)^2 }} \right]^{{n - 2}\over 2}
\end{eqnarray*}
for $x \in {\R}^n \setminus \{ 0 \}$. Hence we have (7.7).\qedwh
{\bf Corollary 7.8.} \ \ {\it Let $u$ be a spherical solution as in proposition 7.6 and $v$ the cylindrical transform of $u$ via (7.2). Then the reflection at $\{ s_o \} \times S^{n - 1}$: $\tilde v (s, \theta) = v (2 s_o - s, \theta)$ is equal to the cylindrical transform of the Kelvin transform of $u$ about the sphere of radius $e^{s_o}$.}\\[0.1in]
{\bf Proof.} \ \ As in the proof of proposition 7.3, we have 
\begin{eqnarray*}
\tilde v (s, \theta) & = &  \left( {{ \lambda e^{2s_o - s}}\over {\lambda^2 + e^{4s_o - 2s} - 2 
e^{2s_o - s} (\theta \cdot \xi) + |\xi|^2}} \right)^{{n - 2}\over 2}\\
& = & \left( {{ \lambda e^{2s_o} e^s}\over {(\lambda^2 + |\xi|^2)e^{2s} +   e^{4s_o} - 2 
e^{2s_o} e^s (\theta \cdot \xi)}} \right)^{{n - 2}\over 2}\\
& = & \left( {{ {{\lambda e^{2s_o}}\over {\lambda^2 + |\xi|^2}} e^s }\over { {{e^{4s_o}}\over {\lambda^2 + |\xi|^2}} + e^{2s}  - {{2 e^{2s_o} e^s (\theta \cdot \xi)}\over {\lambda^2 + |\xi|^2}} }} \right)^{{n - 2}\over 2}\\
& = & \left( {{ \tilde \lambda e^s }\over { {\tilde \lambda}^2 + e^{2s} - 2 e^s (\theta \cdot \tilde \xi) + |\tilde \xi|^2}} \right)^{{n - 2}\over 2}\,, 
\end{eqnarray*}
where
$$
\tilde \lambda = {{e^{2s_o} \lambda }\over {\lambda^2 + |\xi|^2}}  \ \ \ \ {\mbox{and}} \ \ \ \ 
\tilde \xi = {{e^{2s_o} \xi}\over {\lambda^2 + |\xi|^2}}\,.
$$
From proposition 7.6, we can draw the conclusion that $\tilde v$ is the cylindrical transform of a spherical solution $\tilde u$ which is the Kelvin transform of $u$ about the sphere of radius $e^{s_o}$.\qed 
It follows from proposition 7.3 and corollary 7.8 that
$$
v (2 s_o - s, \,\theta) = \tilde v (s, \,\theta) = v (s + \ln (\lambda^2 + |\xi|^2) - 2 s_o, \,\theta)\,.
$$
Thus we have the following symmetric property.\\[0.2in]
{\bf Corollary 7.9.} \ \ {\it Let $u$ be a spherical solution, $v$ the cylindrical transform of $u$, and $\tilde v$ the reflection of $v$ at $\{ s_o \} \times S^{n - 1}$. Then we have
$$
\tilde v (s, \,\theta) = v (2 s_o - s, \,\theta) = v (s + \tilde s, \,\theta) \ \ \ \ for \ \ s \in \R\,, \ \ \theta \in S^{n - 1},
$$
where $\tilde s = \ln (\lambda^2 + |\xi|^2) - 2 s_o\,.$}


\pagebreak

{\bf \Large {\bf 8. \ \ Delaunay-Fowler-type Solutions}}

\vspace{0.3in}

In addition to the spherical solutions, there is a one parameter of smooth positive solutions of equation (6.5) that depend on $s$ only. They pertain to constant mean curvature surfaces and satisfy the ODE 
$$
v'' - \left( {{n -
2}\over 2} \right)^2 v + n (n - 2) v^{{n + 2}\over {n - 2}} = 0 \ \ \ \ {\mbox{in}} \ \ \R\,. \leqno (8.1)
$$
The Hamiltonian energy
$$
H (v, v') = (v')^2 - \left( {{n -
2}\over 2} \right)^2 v^2 + (n - 2)^2 v^{{2n}\over {n - 2}} \leqno (8.2)
$$
is found to be constant along solutions of equation (8.1). It can be seen that when $H (v, v') < 0$, then $v$ cannot be positive everywhere. Let 
$$
D_n = - {2\over n} (n - 2)^2  \left( {{n - 2}\over {4n}} \right)^{{n - 2}\over 2}\,.
$$ 
For any number  
$
H_o \in [D_n , \,0)\,,
$
there exists a unique bounded {\it periodic} smooth positive solution $v$ of equation (8.1) with 
$$
H (v, v') = H_o\,, \ \ \ \ v' (0) = 0 \ \ \ \ {\mbox{and}} \ \ \ \ v'' (0) \le 0\,.
$$ 
There are two extreme cases. When 
$H (v, v') = D_n\,,$ 
$v$ is the constant function given by 
$$
v (t) = \left( {{n - 2}\over {4n}} \right)^{{n - 2}\over 4}\,.
$$
The only positive solution of equation (8.1) with $H (v, v')= 0\,,$ $v' (0) = 0$ and $v'' (0) \le 0$ is given by $v (t) = (2 \cosh t)^{(2-n)/2}\,,$ which is the cylindrical transform of the spherical solution $u_{1, o}$. Smooth positive solutions $v$ of equation (8.1) with 
$H (v, v') \in \left( D_n, \,0\right)$  
are known as {\it Delaunay-Fowler-type solutions}.\bk
Given a Delaunay-Fowler-type solution $v$, let $\varepsilon$  be the minimum of $v (t)\,,$ which is known as the neck-size of the solution. We have
$$
H (v, v') = (n - 2)^2 \varepsilon^{{2n}\over {n - 2}} - \left( {{n -
2}\over 2} \right)^2 \varepsilon^2. \leqno (8.3)
$$  
From the bounds of $H$ we obtain  
$$
0 < \varepsilon < \left( {{n - 2}\over {4n}} \right)^{{n - 2}\over 4}\,.
$$
We may parameterize Delaunay-Fowler-type solutions by $v_{\varepsilon}$, where $\varepsilon$ is the neck-size. Denoted by $T_\varepsilon$ the period of the Delaunay-Fowler-type solution $v_\varepsilon$. Then 
$$
T_\varepsilon \to 0 \ \ \ \ {\mbox{as}} \ \ \ \ \varepsilon \uparrow \left( {{n - 2}\over {4n}} \right)^{{n - 2}\over 4} \ \ \ \ {\mbox{and}} \ \ \ \  T_\varepsilon \to \infty \ \ \ \ {\mbox{as}} \ \ \ \ \varepsilon \downarrow 0
$$
(see \cite{Mazzeo-Pacard_2}). Moreover, $v (t)$ converges to $(2 \cosh t)^{(2-n)/2}$ in compact subsets of ${\R}^n$ as $\varepsilon \downarrow 0$ (\cite{K-M-P-S}, cf. \cite{Leung.5}).\bk
By a result of Caffarelli-Gidas-Spruck \cite{Caffarelli-Gidas-Spruck}, if $v$ is a smooth positive solution of equation (6.5) in $\R \times S^{n - 1}$,\, and if $v$ has a ``removable" singularity at $t = - \infty$ (in the sense that the corresponding function has a removable singularity at $x = 0$), then $v$ is the cylindrical transform of one of the spherical solutions. Moreover, if  $v$ has a non-removable singularity at $t = -\infty$, then it is one of the Delaunay-Fowler-type solutions. When we consider the equation outside a compact set of ${\R}^n$ and $K$ is equal to one for $t \gg 1$, we have more variety by deforming the Delaunay-Fowler-type solutions.\bk
Let $u$ be a smooth positive solution of (4.1) and $v$ be defined by the equation
$$
u (x) = |x|^{{2-n}\over 2} v \left(- \ln |x|\,, \,{x\over {|x|}} \right) \mfor x \in {\R}^n \setminus \{ 0 \}\,. \leqno (8.4)
$$
Applying the Kelvin transform about the sphere of radius $a > 0$ we obtain
$$
|x|^{{2-n}\over 2}  v \left( \ln {{|x|}\over a^2}\,, \,{x\over {|x|}}\right) \mfor x \in {\R}^n \setminus \{ 0 \}\,. \leqno (8.5)
$$
By a translation $x \mapsto x + \xi$, we have
$$
|x - \xi|^{{2-n}\over 2} v \left( \ln |x - \xi| - \ln a^2\,, \,{{x - \xi}\over {|x - \xi|}} \right) \mfor x \in {\R}^n \setminus \{ \xi \}\,, \leqno (8.6)
$$
where $\xi \in {\R}^n.$ Applying the Kelvin transform about the sphere of radius $a$ again we obtain
\begin{eqnarray*}
(8.7) & \, &  u_{\xi, a} (x)\\
 & := & |x|^{2 - n} a^{n - 2}  \,\bigg\vert {{a^2x}\over {|x|^2}} - \xi \bigg\vert^{{2-n}\over 2} \,v \left( \ln  \bigg\vert {{a^2x}\over {|x|^2}}- \xi \bigg\vert - \ln a^2\,, \,\left( {{a^2x}\over {|x|^2}} - \xi \right)/\bigg\vert{{a^2x}\over {|x|^2}}- \xi \bigg\vert \right)\\
& = & |x|^{{2 - n}\over 2} \,\bigg\vert \theta - {{\xi}\over a^2} |x| \bigg\vert^{{2-n}\over 2} \,
v \left( - \ln |x| + \ln \bigg\vert \theta - {{\xi}\over a^2} |x| \bigg\vert\,, \,
\left( \theta - {{\xi}\over a^2} |x| \right) /\bigg\vert  \theta - {{\xi}\over a^2} |x|  \bigg\vert
\right)  
\end{eqnarray*}
for $x/|x|^2 \not= \xi/a^2\,,$ where $\theta = x/|x|$ for $x \not= 0$. The function $u_{\xi, a}$ satisfies the equation 
$$
\Delta u_{\xi, a} + K_{\xi, a} \, u_{\xi, a}^{{n + 2}\over {n - 2}} = 0\,, \leqno (8.8)
$$
where
$$
K_{\xi, a} (x) = K \left(  a^2\left( {{a^2 x}\over {|x|^2}} - \xi \right)/\bigg\vert  {{a^2 x}\over {|x|^2}} - \xi \bigg\vert^2  \right) 
= K \left( \left( {{x}\over {|x|^2}} - {\xi \over {a^2}} \right)/\bigg\vert  {{x}\over {|x|^2}} - {\xi \over {a^2}} \bigg\vert^2  \right)\,. \leqno (8.9)
$$
It turns out that the solution $u_{\xi, a}$ and the equation depend on $\xi/a^2$ only. We observe that $K_{\xi, a}$ is equal to one outside a compact subset in ${\R}^n$. In case $v$ is one of the Delaunay-Fowler-type solutions, the above deformation is first expounded in \cite{K-M-P-S}. In this case we may let $a = 1$, $|x| = e^{-t}$ and 
$$
v_{\varepsilon, \xi} (t, \theta) = | \theta -  \xi e^{-t} |^{{2-n}\over 2} \,
v_{\varepsilon} ( t  + \ln | \theta - \xi e^{-t} |)
$$
(cf. the last line in (8.7)). For $T \in {\R}$, define the deformed Delaunay-Fowler-type solution by 
$$
v_{\varepsilon, \xi, T} (t, \theta) = v_{\varepsilon, \xi} (t + T, \theta)\,, \leqno (8.10)
$$
and set  
$$
u_{\varepsilon, \xi, T} (x) = |x|^{{2-n}\over 2} 
v_{\varepsilon, \xi, T} (-\ln |x|, \theta) \mfor x \in {\R} \setminus \{ 0 \}   \leqno (8.11) 
$$
(cf. (8.7)).

\vspace*{0.4in}


{\bf \Large {\bf 9. \ \ Pohozaev Invariants}}

\vspace{0.3in}

We further explore the conformal structure of equation (4.1). 
For a compact Riemannian $n$-manifold $(N, g)$ with smooth boundary $\partial N\,,$ $n \ge 3$,  let $R$ be the scalar curvature of $(N, g)$, $X$ a conformal Killing vector field on $N$. In \cite{Schoen.2} (cf. also \cite{Pacard-Riviere}), Schoen obtains the identity  
$$
{{(n - 2)^2}\over {8n (n - 1)}} \int_N ({\cal L}_X R) \, dv_g = {{n - 2}\over {4(n - 1)}} \int_{\partial N} T (X, \nu) \, d \sigma\,, \leqno (9.1)
$$
where $T (\cdot\,, \cdot)$ is the trace-free Ricci tensor for the metric $g$, $\nu$ the unit outward normal on $\partial N$, and ${\cal L}_X$ the Lie derivative. (The constants are left in equation (9.1) for later simplification.)\bk 
We apply the identity to the equation 
$$
\Delta u + K u^{{n + 2}\over {n - 2}} = 0 \ \ \ \ {\mbox{in}} \ \ {\R}^n, \leqno (9.2)
$$
which is the same as equation (4.1), except the constant $n(n - 2)$. The  scalar curvature of the conformal metric $u^{4/(n - 2)} \delta_{ij}$ is equal to $4K (n - 1)/(n - 2)$. Here $\delta_{ij}$ is Euclidean metric. Take $N = B_o (r)$. The left hand side of (9.1) becomes
$$
{{n - 2}\over {2n}} \int_{B_o (r)} ({\cal L}_X K) \, u^{{2n}\over {n - 2}} dx\,, \leqno (9.3)
$$
and the right hand side is given by (see \cite{Schoen.2})
\begin{eqnarray*}
(9.4) \ \ \ \ \ \ \ \ & \ & {{n - 2}\over {4 (n - 1)}} \int_{S_r} T (X, \nu) \, d \sigma\\
& = & \int_{S_r} |x|^{n - 2} \left\{ {n\over {n - 1}} (X \cdot \btd u) ({\cal D} \cdot \btd u) - {{n - 2}\over {n - 1}} u \,(D du) (X, {\cal D}) \right. \ \ \ \ \ \ \ \ \ \ \ \ \ \ \ 
\\
& \ & \ \ \ - \left. \left( {1\over {n - 1}} |\btd u|^2 + {{(n - 2)^2}\over {4 (n - 1)}} u^{{2n}\over {n - 2}} \right) (X \cdot {\cal D})\right\} \, dS\,,
\end{eqnarray*}
where ${\cal D} = \sum x_i \partial_{x_i}\,,$ and $D d u$ is the Hessian of $u$. Take $X = {\cal D}$, it follows from a direct calculation that 
\begin{eqnarray*}
(9.5) \ \ \ \ \ \ \ \ & \ & {{n - 2}\over {2n}} 
\int_{B_o (r)} x
\cdot
\btd K (x) \,u^{{2n}\over {n - 2}} (x)\, dx\\
& = & \int_{S_r} \left[ r \left( {{\partial u}\over {\partial r}} \right)^2
 - {r\over 2} |\btd u|^2 +  {{n - 2}\over {2n}}\, r  K u^{{2n}\over {n -
2}}  + {{n - 2}\over 2} \,u \,{{\partial
u}\over {\partial r}} \right] \,dS\,, \ \ \ \ \ \ \ \ \ \ \ \ \ \ \ 
\end{eqnarray*}
which is known as the {\it radial Pohozaev identity}. Denote the right hand side of (9.5) by $P (u, r)$. If we let $X = \partial_{x_i}$, we obtain 
\begin{eqnarray*}
(9.6) \ \ \ \ \ \ \ \ \ \ & \ & {{n - 2}\over {2n}} \int_{B_o (r)} {{\partial K}\over {\partial x_i}} (x) \, u^{{2n}\over {n - 2}} (x) \, dx\\
 & = & \int_{S_r} \left[ {{\partial u}\over {\partial x_i}} {{\partial u}\over {\partial r}} - {1\over 2} | \btd u|^2 \nu_i + {{n - 2}\over {2n}} K u^{{2n}\over {n - 2}} 
\nu_i \right] \, dS\,,  \ \ \ \ \ \  \ \ \ \ \ \ \ \ \ \ \ \ \ \ \ \ \ \ \ \ \ \  
\end{eqnarray*}
where $\nu_i$ is the $i$-th component of the unit outward normal $\nu$ on $S_r$. (9.6) is called  the {\it translational Pohozaev identity}. Denote the right hand side of (9.6) by $P_i (u, r)\,.$\bk
Assume that $u$ has slow decay, that is,
$$
u (x) \le C |x|^{(2 - n)/2}
$$
for large $|x|$. From lemma 4.4 we have the gradient estimate
$$
|\btd u (x)| \le C' |x|^{-n/2}
$$
for large $|x|$. Here $C$ and $C'$ are positive constants. The integrand in the right hand side of (9.6) is of order $O(|x|^{-n})$. It follows that 
$\displaystyle{\lim_{r\to \infty} P_i (u, r) = 0}$ and hence 
$$
\lim_{r \to \infty} \int_{B_o (r)} {{\partial K}\over {\partial x_i}} (x) \, u^{{2n}\over {n - 2}} (x) \, dx = 0 \mfor i = 1, \,2,..., n\,.
$$
\hspace*{0.5in}However, under the Kelvin transform, (9.6) could yield a non-zero invariant. Let 
$$
\tilde u (x) = |x|^{2 - n} \,u (x/|x|^2) \mfor x \not= 0
$$
be the Kelvin transform of $u$. $\tilde u$ satisfies the equation
$$
\Delta \tilde u + K (x/|x|^2) \,{\tilde u}^{{n + 2}\over {n - 2}} = 0 \ \ \ \ {\mbox{in}} \ \ x \in {\R}^n \setminus \{ 0 \}\,.
$$
If $u$ has slow decay, then 
$$
\tilde u (x) \le C |x|^{(2 - n)/2} \ \ \ \ {\mbox{and}} \ \ \ \ |\btd \tilde u (x)| \le C' |x|^{-n/2}
$$
for small $|x|$ (cf. lemma 4.4). As in (9.6) we have
$$
{{n - 2}\over {2n}} \int_{B_o (r) \setminus B_o (s)} {{\partial \tilde K}\over {\partial x_i}} (x) \, {\tilde u}^{{2n}\over {n - 2}} (x) \, dx = P_i (\tilde u, r) - P_i (\tilde u, s) \leqno (9.7)
$$
for $ r > s > 0\,,$ where $\tilde K (x) = K (x/|x|^2)$. In particular, if 
$$
\bigg\vert {{\partial \tilde K}\over {\partial x_i}} (x) \bigg\vert  \le C\, |x|^\alpha
$$
for small $|x|$ and some positive constants $C$ and $\alpha$, then 
$$
\int_{B_o (r)} {{\partial \tilde K}\over {\partial x_i}} (x) \, {\tilde u}^{{2n}\over {n - 2}} (x) \, dx
$$
is finite. It follows that the limit  
$$
\lim_{s \to 0^+} P_i (\tilde u, s) =  P_i (\tilde u, r) - {{n - 2}\over {2n}} \int_{B_o (r)} {{\partial \tilde K}\over {\partial x_i}} (x) \, {\tilde u}^{{2n}\over {n - 2}} (x) \, dx \leqno (9.8)
$$
exists and is known as the {\it $i$-th translational Pohozaev invariant}, which is denoted by $P_i (\tilde u)$.\bk
For the radial Pohozaev invariant,  the expression $P (u, r)$ is invariant under the Kelvin transform, that is,  
$$
P (\tilde u, s) = P (u, r) \mfor r = 1/s\,.
$$
Hence 
$$
 \lim_{r \to + \infty} P (u, r) = \lim_{s \to 0^+}  P (\tilde u, s)\,. \leqno (9.9)
$$
The limit in (9.9), if it exists, is called the {\it radial Pohozaev invariant} and is denoted by $P (\tilde u)$. It is found that for $u_{\varepsilon, a, T}$ in (8.11), 
$$
P (u_{\varepsilon, a, T})  = {{(n - 2)^2 }\over 8}\, \omega_{n -1} \left(4 \varepsilon^{{2n}\over {n - 2}} - \varepsilon^2 \right) = {\omega_{n -1} \over 2} H (\varepsilon)\,, \leqno (9.10)
$$
and 
$$
P_i (u_{\varepsilon, a, T}) = {2n\over {n - 2}} \omega_{n - 1} H (\varepsilon)\, \xi_i = {{4n}\over {n - 2}} P (\tilde u) \,\xi_i \mfor i = 1,..., n \leqno (9.11)
$$
(cf. \cite{K-M-P-S}; here $K (x) = n (n - 2)$ for small $|x|$).

\vspace*{0.4in}


{\bf \Large {\bf 10. \ \ Questions on upper and lower bounds,  \hfill 
\vspace*{0.05in}
\hspace*{0.55in} Existence and Uniqueness}}

\vspace{0.3in}

As a first step, one would like to find conditions on $K$ so as to ensure slow decay (i) of a smooth positive solution $u$ of equation (4.1). This appears to be essential 
for a better understanding on the asymptotic behaviour of positive solutions $u$. When $K$ is only required to be bounded between two positive numbers, there may have solutions which are unbounded from above (see \cite{Taliaferro.1} and \cite{Leung.4}, cf. also \cite{Taliaferro.2}). The works of Chen and Lin \cite{Chen-Lin.1} \cite{Chen-Lin.2} \cite{Chen-Lin.3} \cite{Chen-Lin.4}  
exhibit conditions on $K$ for slow decay. It is conjectured by Taliaferro in \cite{Taliaferro.1} that if $K$ satisfies condition (4.2) with $b/a < 2^{2/(n - 2)}$,\, then $u$ should have slow decay.\bk
Once a solution $u$ of equation (4.1) has slow decay, a natural next step is to ask if it has lower bound (iii) as well. It turns out that this is not always the case (see \cite{Chen-Lin.4}  and  \cite{Leung.5}), and such $u$ are called {\it exotic solutions}. In \cite{Lin.1} Lin states a conjecture on the existence of exotic solutions based on a nondegenerate condition related to $K$.\bk
Another concern is uniqueness. Assume that $K \equiv 1$ outside a compact subset of ${\R}^n.$  Let $u_1$ and $u_2$ be smooth positive solutions of equation (4.1). It is shown in \cite{K-M-P-S} that there exist deformed Delaunay-Fowler-type solutions $v_{\varepsilon_1, \xi_1, T_1}$ and $v_{\varepsilon_2, \xi_2, T_2}$ (cf. (8.10) and (8.11)) such that 
$$
u_i (x) = |x|^{{2-n}\over 2} 
v_{\varepsilon_i, \xi^i, T_i} \,(\ln |x|,\, \theta)  \,+ \,O (|x|^{-\alpha}) \mfor |x| \gg 1\,, \ \ i = 1, \ 2\,,
$$
where $\xi^1\,, \ \xi^2 \in {\R}^n$, and $\alpha$ is a positive constant. The parameters  $\varepsilon_i$ and $\xi_i$ can be determined by the Pohozaev invariants of $u_i$, $i = 1, \ 2,$ respectively. (For the translational Pohozaev invariants, we take the Kelvin transforms of $u_1$ and $u_2$, cf. section  9.) It is asked (cf. \cite{K-M-P-S}) when $\varepsilon_1 = \varepsilon_2,$ $\xi^1 = \xi^2$ and $T_1 =T_2$, do we have $u_1 = u_2$?


\pagebreak

\end{document}